\documentclass[12pt]{article}
\usepackage{amsmath}
\usepackage{amsfonts}
\usepackage{amssymb}
\usepackage{amsthm}
\usepackage{stackrel}

\usepackage{color}

\usepackage{figsize}

\usepackage{graphics}

\newtheorem{thm}{Theorem}[section]

\newcommand{\R}{{\rm I}\kern-0.18em{\rm R}}
\newcommand{\1}{{\rm 1}\kern-0.25em{\rm I}}
\newcommand{\E}{{\rm I}\kern-0.18em{\rm E}}
\newcommand{\p}{{\rm I}\kern-0.18em{\rm P}}
\makeatletter

\title{Big Outliers Versus Heavy Tails: what to use?}
\author{Lev B. Klebanov\\
Department of Probability and Mathematical Statistics,\\ Charles University}
\date{}
\begin{document}
\maketitle

\begin{abstract}
A possibility to give strong mathematical definitions of outliers and heavy tailed distributions or their modification is discussed. Some alternatives for the notion of tail index are proposed.

Key words: outliers, heavy tails, tail index.
\end{abstract}

\section{Intuitive approach or mathematical definition?}\label{s1}
\setcounter{equation}{0}

Professor Jerzy Neyman in his talk ``Current Problems of Mathematical Statistics" \cite{JN} wrote: ``In general, the present stage of development of mathematical statistics may be compared with that of analysis in the epoch of Weierstrass." Although we have many new mathematically correct results in Statistics, the situation seems to be similar now. There are some ``intuition-made" definitions of objects that have no precise sense in Statistics. The use of such definitions seems sometimes very strange. Here I would like to discuss two of such objects: {\bf outliers} and  {\bf heavy tails}. 

\vspace{0.2cm}
Let us start with {\bf heavy tails}. At the first glance, the notion seems to be clear and nice. Really, if $X$ is a random variable (r.v.) then its tail is defined by the relation
\begin{equation}\label{eq1}
T(x)=T_X(x) = \p\{|X|>x\},\; \; \; x>0.
\end{equation}
Obviously, the definition of the tail $T(x)$ is absolutely correct. 

However, what does it mean that the tail is heavy? One of used definitions is the following.
We say r.v. $X$ has heavy (power) tail with parameters $\alpha>0$ and $\lambda>0$ if there exists the limit
\begin{equation}\label{eq2}
\lim_{x \to \infty}T(x) x^{\alpha} =\lambda .
\end{equation} 
Let us look at (\ref{eq2}) more attentively. If we have two different r.v.s $X$ and $Y$ such that
$T_X(x)=T_Y(x)$ for all $x>A$, where $A$ is a positive number, then all parameters $\alpha$ and $\lambda$ in (\ref{eq2}) are the same for both $T_X$ and $T_Y$ that is both $X$ and $Y$ have heavy tail with parameters $\alpha$ and $\lambda$. We say that r.v.s are equivalent if their tails are identical in a neighborhood of infinity. Then we may talk about classes of equivalence for all r.v.s. All r.v.s from each equivalence class have (or do not have) heavy tail with the same parameters. 

\vspace{0.2cm}
What does it mean from statistical point of view? It means that (for non-parametric situation) {\bf we can never estimate the parameters $\alpha$ and/or $\lambda$}. Really, for each finite set $x_1, \ldots , x_n$ of observations on r.v. $X$ we can never say what will be the behavior of $T(x)$ for $x > \max{|x_1|, \ldots , |x_n|}$.

To have a possibility of such estimation we need either to restrict ourselves with a small class of r.v.s under consideration, or modify the notion of heavy tail. Of course, we need mathematically correct definition which is suitable for statistical study. However, we shall go back to this problem a little bit later.

\vspace{0.2cm}
Let us consider a notion of {\bf outliers} now. It is one of the most strange notions from my view. Wikipedia, the free encyclopedia defines outliers in the following way: ``In statistics, an outlier is an observation point that is distant from other observations. An outlier may be due to variability in the measurement or it may indicate experimental error; the latter are sometimes excluded from the data." I think, some points from this definition need essential clarification. Really, let us consider the following graphs.

\begin{figure}[htp]
	\centering
	\hfil
	\includegraphics[scale=0.7]{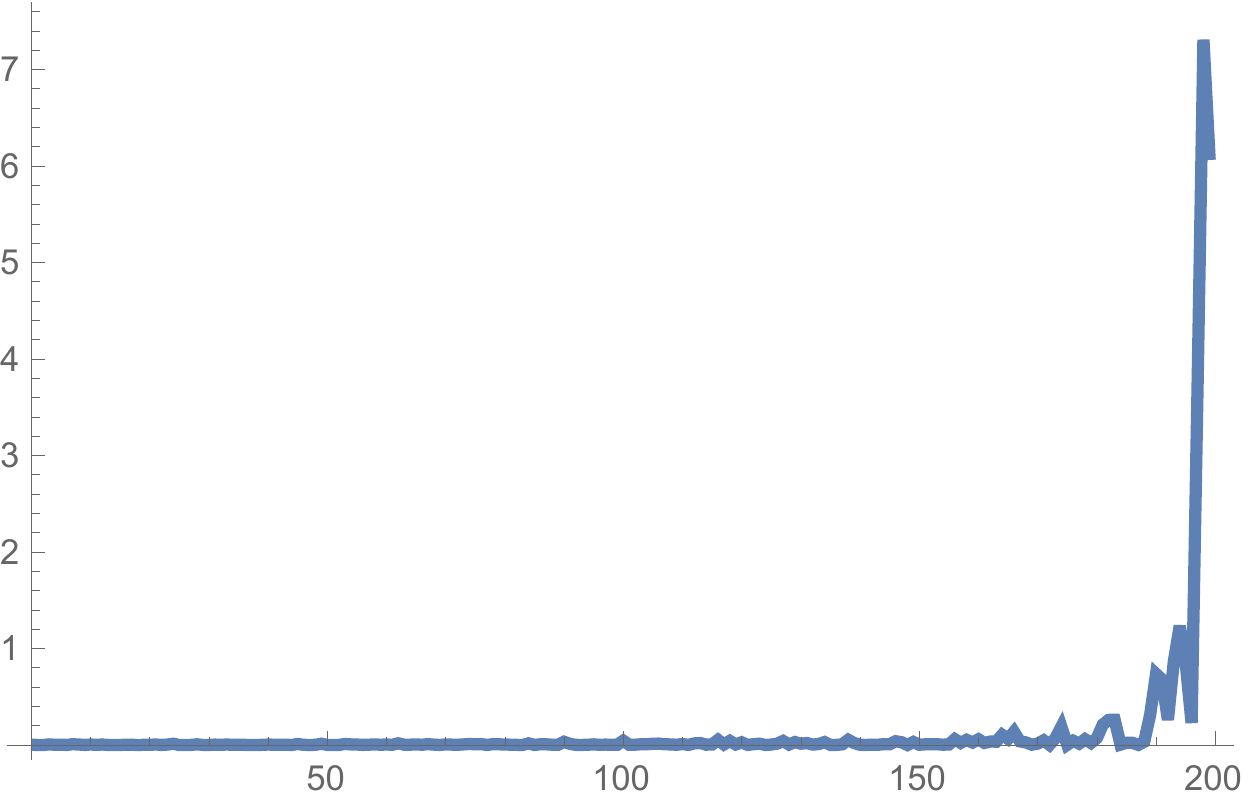}
	\caption{Distances between ordered statistics for the sample of volume 200 from Pareto distribution (0,2)}\label{fig1}
\end{figure}

On Figure \ref{fig1} the distance between $|X|_{n,n}$ and $|X|_{n-1,n}$ is greater that ``typical" distance between order statistics in 40-50 times. So, it seems (intuitively) we have outliers here. Of course, it is in intuitive agreement with the fact the sample was taken from Pareto distribution.

	\begin{figure}[htp]
		\centering
		\hfil
		\includegraphics[scale=0.7]{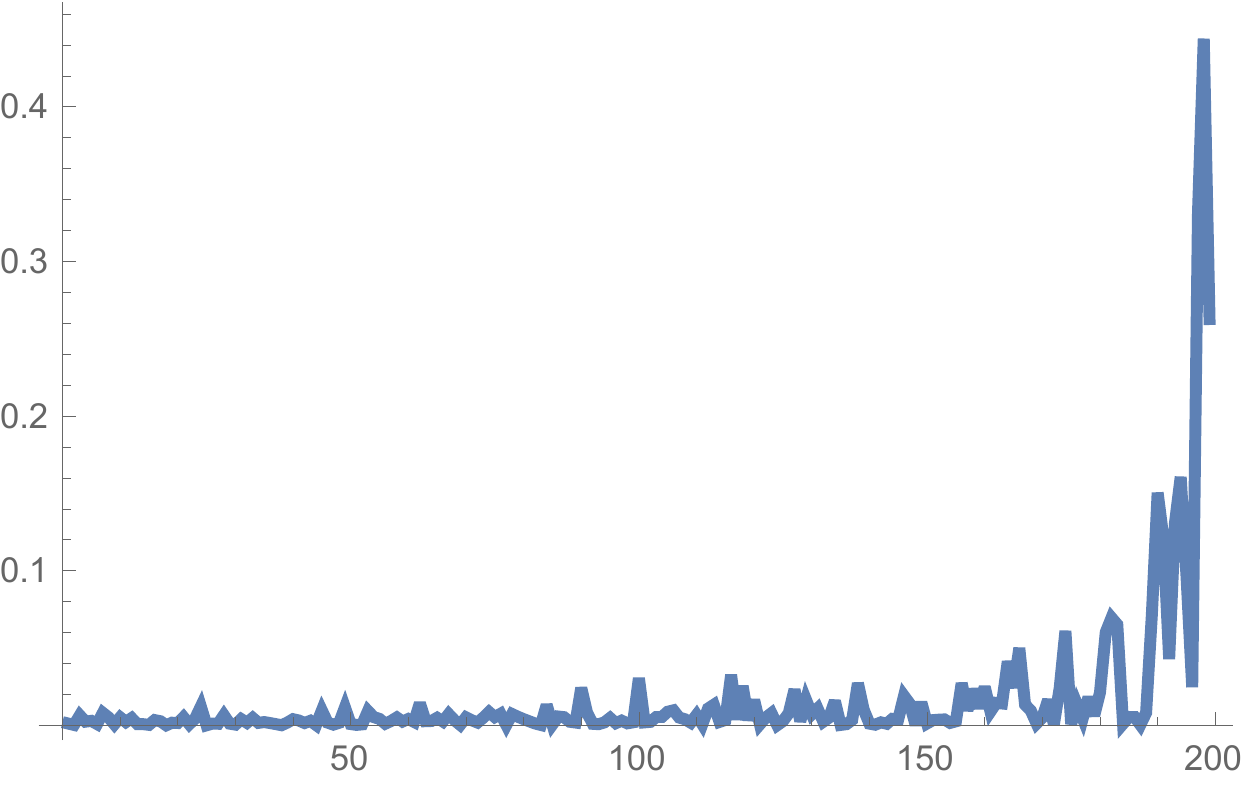}
		\caption{Distances between logs of ordered statistics for the sample of volume 200 from Pareto distribution (0,2)}\label{fig2}
	\end{figure}
On Figure \ref{fig2} the distance between $|X|_{n,n}$ and $|X|_{n-1,n}$ is greater that ``typical" distance between order statistics in 30-35 times. It is smaller that for previous case. However, without comparing this with Figure \ref{fig1} we cannot say 30-35 times is not large enough. Intuitively, we have outliers again. However, the sample now is from exponential distribution, which is not heavy-tailed.

	\begin{figure}[htp]
		\centering
		\hfil
		\includegraphics[scale=0.7]{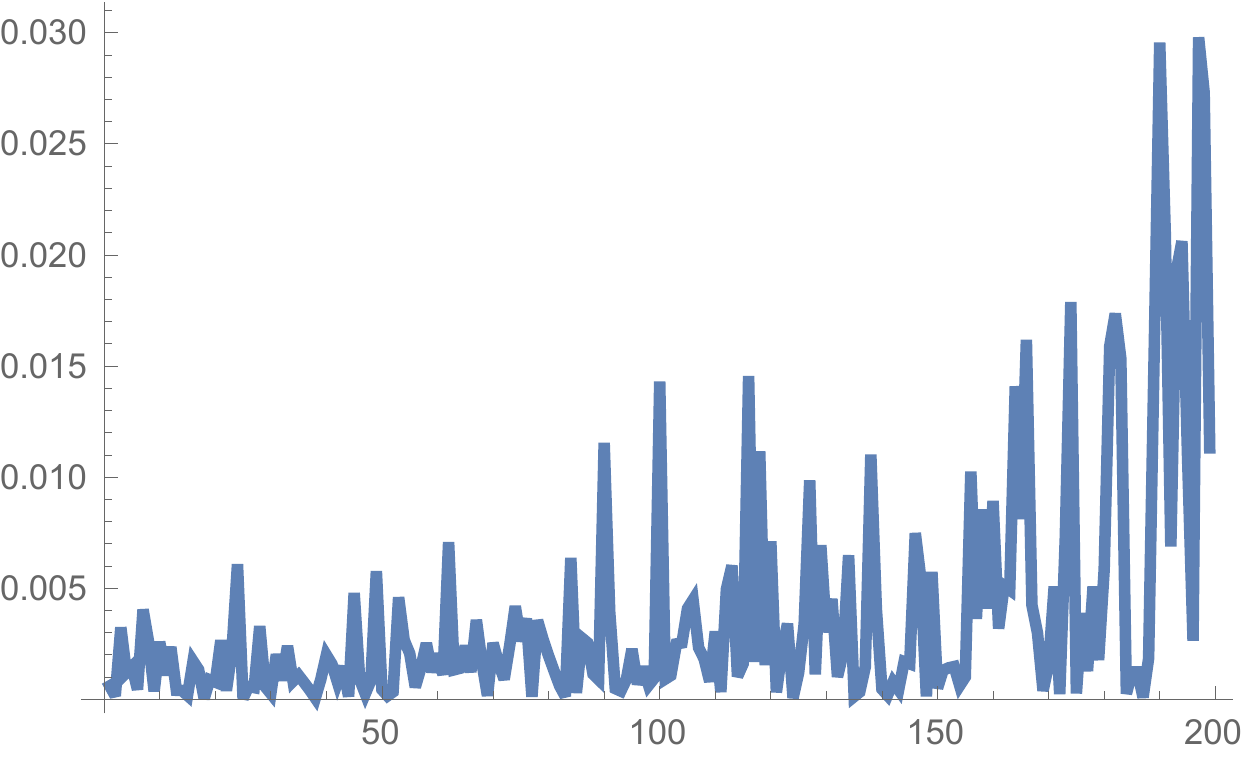}
		\caption{Distances between ArcTan's of ordered statistics for the sample of volume 200 from Pareto distribution (0,2)}\label{fig3}
	\end{figure}	

On Figure \ref{fig3} we have sample from a distribution with compact support. It is neither heavy-tailed nor high variability (in terms of large standard deviation). However, we see that the difference between $|X|_{n,n}$ and $|X|_{n-1,n}$ is greater that ``typical" distance between order statistics 10-15 times. How can we see it is not enough to say about outliers?

More generally, {\bf how is it possible to discuss the presence of outliers, if we always can transform arbitrary r.v.s in corresponding set of bounded random variables without loss of statistical information?}

\vspace{0.2cm} The answer is simple. Usually, statisticians study a scheme in which r.v.s are  generated. If we like to transform r.v.s it is necessarily to change the scheme in corresponding way, which may be not too easy. For example, if we study sums of r.v.s
\[ S_n= X_1 + \ldots +X_n\]
the transformation from $X_j$ to, say, $\arctan X_j$ will change summation of $X_j$ to an unclear operation.

\vspace{0.2cm} This leads us to an idea that {\bf the notion of outliers has to be considered not by itself but in connection with underlying scheme.} If so, we must study different schemes, leading to some sets of r.v.s, especially to that with heavy-tailed distributions. 

\subsection{Characterizations of r.v.s}
I think that such schemes and corresponding distributions of r.v.s are natural products by Characterization of Probability distributions. The aim of Characterizations is to describe all distributions of random variables possessing a desirable property, which may be taking as a base of probabilistic and/or statistical model.

\vspace{0.2cm}
Let us start with an example leading to Polya Theorem \cite{Pol}. Suppose that we have a gas whose molecules are chaotically moving, and the space is isotropic and homogeneous. Denote by $X_1$ and $X_2$ projections of the velocity of a molecule on the axis in $(x,y)$ plain. In view of space property we have the following properties: a) $X_1$ and $X_2$ are independent random variables; b) $X_1 \stackrel{d}{=}X_2$. After rotation of the coordinate system counter clock wise on the angle $\pi/4$ we obtain, that a projection on new coordinate axes has to be identically distributed with the old one. That is, $X_1 \stackrel{d}{=}(X_1+X_2)/\sqrt{2}$. {\it  Polya Theorem says that in this situation $X_1$ has normal (or degenerate) distribution with zero mean}.

From Polya Theorem we obtain Maxwell distribution for velocities of gas molecules basing on two natural properties of the space as isotropy and homogeneity only. {\bf Are there any models leading in a natural way to heavy-tailed distributions?}

\vspace{0.2cm}
Let us show, that strictly stable distributions may be also described by a clear physical property. {\bf Let us explain this by an example taken from mobile telephoning}: Suppose that we have a base station. And suppose that there is a Poisson ensemble of points (Poisson field), the locations of mobile phones. Each phone produces a random signal $Y_k$. It is known that the signal depression is in inverse proportion with a power of the distance $G_k$ from the phone to base station. Therefore, the cumulative signal coming to base station can be represented as $X =Y_1/G{_1}^a+...+Y_n/G_{n}^a+...$. This is LePage series \cite{LeP}, and it converges to a strictly stable distribution with index $\alpha=1/a$.  Obviously, we may change the base station and mobile phones by electric charges, or by physical masses. In any such case we obtain stable non-normal distribution of the resulting forth.  

Heaviness of tail for strictly stable distribution is defined by the index of stability $\alpha$, which may be expressed through signal depression.	The last is a physical characteristic which can be estimated directly (not through observations of $X$).

\vspace{0.2cm}
It is clear, that for this scheme there will be many observations on $X$, which seem to be ``far" from each other. But is it natural to call them ``outliers"? Do they indicate experimental errors? Definitely, the answer to the last question is negative. On the other hand, variability of the measurements here is high, but natural. I think, we have no reasons to consider such observations as something special, to what one need pay additional attention. Of course, we may not ignore such observations.

\subsection{Toy-model of capital distribution}
	In physics, under toy-model usually understand a model, which does not give complete description of a phenomena, but is rather simple and provides explanation of essential part of the phenomena.
	
	Let us try to construct a toy-model for capital distribution (see \cite{KMR}). 
	Assume that there is an output (business) in which we invest a unit of the
	capital at the initial moment $t=0$. at the moment $t=1$ we get a sum of
	capital $X_1$ (the nature of the r.v. $X_1$ depends on the nature  of
	the output and that of the market). If the whole sum of capital remains in the
	business, then to the moment $t=2$ the sum of capital becomes $X_1\cdot X_2$,
	where r.v. $X_2$ is independent of $X_1$ and has the same distribution as
	$X_1$ (provided that conditions of the output and of the market are
	invariable). Using the same arguments further on, we find that to the moment
	$t=n$ the sum of capital equals to $\prod_{j=1}^{n}X_j$, and also r.v.s
	$X_1, \ldots ,X_n$ are i.i.d. 

	From the economical sense it is clear that
	$X_j > 0$, $j=1, \ldots ,n$. Now assume that there can happen a change
	of output or of the market conditions which makes further investment of
	capital in the business impossible. We assume that the time till the appearance of the unfavorable
	event is random variable $\nu_p$, $p=1/\E \nu_p$.
	The sum of capital to the moment of this event equals to
	$\prod_{j=1}^{\nu_p}X_j$. And the mean time to the appearance of the unfavorable
	event is $\E \nu_p = 1/p$. Therefore ``mean annual sum of capital" is
	\[ Z_p = \Bigl(\prod_{j=1}^{\nu_p}X_j\Bigr)^p. \]
	The smaller is the value of $p>0$ the rarely is the unfavorable event. If $p$ is small enough, we may approximate the distribution of $Z_p$ by its limit distribution for $p \to 0$. To find this distribution it is possible to pass from $X_j$ to $Y_j = \log X_j$, and change the product by a sum of random number $\nu_p$ of random variables $Y_j$.

	If probability generating functions of $\nu_p$ generate a commutative semigroup, the limit distribution of the sum will coincide with $\nu$-stable or with $\nu$-degenerate distribution. 
	
	1. The most simplest case is that of geometric distribution of $\nu_p$. In this situation, the probability of unfavorable event is the same for each time moment $t=k$. If there exists positive first moment of $Y_j= \log X_j$, then the limit distribution of random sum coincides with $\nu$-degenerate distribution, and is Exponential distribution. This means, that limit distribution of $Z_p$ is Pareto distribution $F(x) = 1-x^{-1/\gamma}$ for $x>1$, and $F(x)=0$ for $x \leq 1$. Here $\gamma = \E \log X_1 >0$. This distribution has power tail. For $\gamma \geq 1$ this distribution has infinite mean.  
	Pareto distribution was introduced by Wilfredo Pareto to describe the capital distribution, but he used empirical study only, and had no toy-model. About hundred years ago this distribution gave a very good agreement with observed facts. Nowadays, we need a small modification of the distribution. Let us mention that our toy-model shows, that such distribution of capitals may be explained just by random effects. This is an essential argument against Elite Theory, because the definition of elite becomes not clear. 

\begin{figure}[htp]
	\centering
	\hfil
	\includegraphics[scale=0.7]{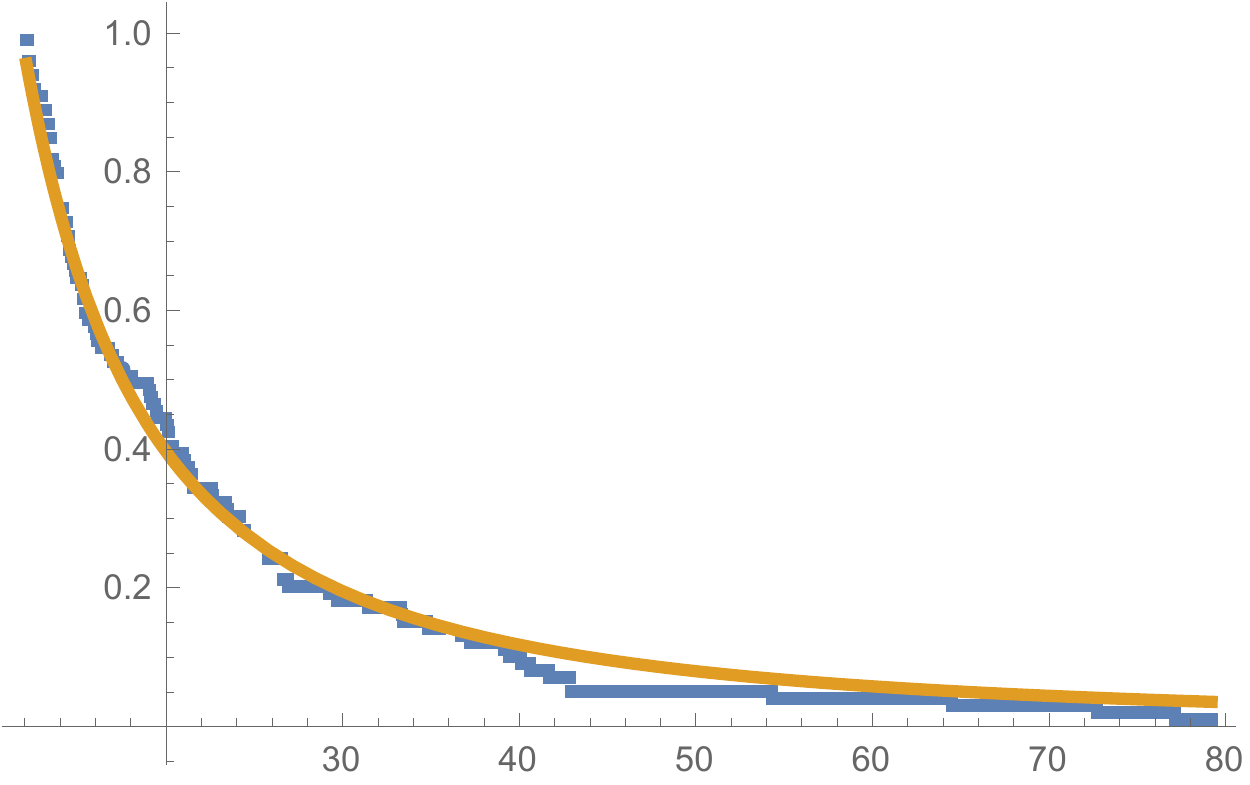}
	\caption{Plot of Pareto distribution function versus empirical distribution of the capital of highest 100 billionaires. Forbes dataset.}\label{fig4}
\end{figure}

\begin{figure}[htp]
	\centering
	\hfil
	\includegraphics[scale=0.7]{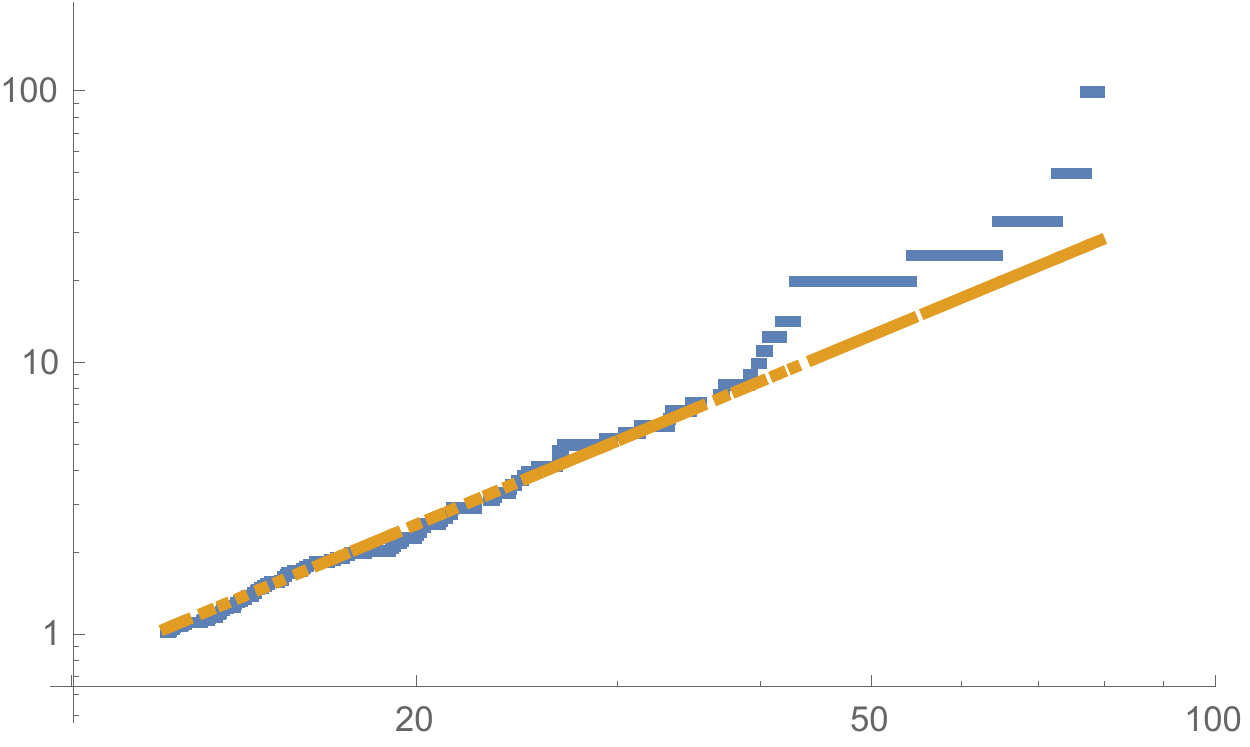}
	\caption{Log-Log-plot of Pareto distribution function versus empirical distribution of the capital of highest 100 billionaires. Forbes dataset.}\label{fig5}
\end{figure}

The situation in the model of capital distribution is, in some sense, similar to that in mobile telephoning model. Namely, statistician will observe large distances between order statistics, but he/she will have no reasons to consider corresponding observations as something special. To estimate the parameter of tail heaviness it is enough to construct an estimator of $\E \log X$. Such estimator is $(1/n)\sum_{j=1}^n \log X_j$. Very important fact is that r.v. $X_j$ may be just bounded while the limit Pareto distribution has power (heavy) tail.

\vspace{0.2cm} Let us note that a very similar model may be obtained through change of the product of random variables $X_j$ by their random number minimum. Again, the r.v. $X_j$ may be bounded, but the limit distribution has heavy tail.

\vspace{0.2cm} It is also of essential interest that such situation is impossible for sums of r.v.s. For limit distribution to have heavy tail it is necessary the summands must have heavy tails too.

Remarkable that for the cases of random products, random minimums and random sums we have the same equation and the same solution for {\it different} transforms of distribution function. They are Mellin transform, survival function and characteristic function correspondingly. I think, it is essential for teaching both Probability and Statistics. The idea to use different transformation of distribution function to get characterization and/or limit theorem is very fruitful, and attempts to omit teaching of, say, characteristic function seems to be just bad simplification of the course of Probability.

\vspace{0.2cm}
{\bf Let us went back to the notion of outliers.} In the definition given above we are talking on some observations ``distant" from other points. What is the ``unit of measurement" for such distance? There are attempts to measure the distance from an observation to their mean value in term of sample variance. 

Suppose that $X_1,X_2, \ldots ,X_n$ is a sequence of i.i.d. r.v.s. Denote by
\[ \bar{x}_n=\frac{1}{n}\sum_{j=1}^{n}X_j, \;\; s^2_n=\frac{1}{n}\sum_{j=1}^{n}(X_j-\bar{x})^2\]
their empirical mean and empirical variance correspondingly. Let $k>0$ be a fixed number. 
Namely, let us estimate the following probability 
\begin{equation}\label{eq3}
p_n=\p \{|X-\bar{x}_n|/s_n>k\},
\end{equation}
It is recommended to say that the distribution of $X$ produces many outliers if the probability (\ref{eq3}) is high (say, higher than for normal distribution). 

The observations $X_j$ for which the inequality $|X_j-\bar{x}_n|/s_n>k$ holds are called outliers. Unfortunately, this approach appears to be not connected to heavy-tailed distributions
(see \cite{KV}).
\vspace{0.2cm}

\begin{thm}\label{th1} \ref{th1}. Suppose that $X_1,X_2, \ldots ,X_n$ is a sequence of i.i.d. r.v.s belonging to a domain of attraction of strictly stable random variable with index of stability $\alpha \in (0,2)$. Then
	\begin{equation}
	\label{eq4}
	\lim_{n \to \infty}p_n =0.
	\end{equation}
\end{thm} 
From this Theorem it follows that (for sufficiently large $n$) many heavy-tailed distributions will not produce any outliers. This is in contradiction with our wish to have outliers for distributions with high variance.
By the way, the word {\bf variability} is not defined precisely, too. It shows, that high variability may denote something different than high standard deviation. Namely, one can observe outliers when the density posses a high peak. 

\section{How to obtain more outliers?}\label{s2}
\setcounter{equation}{0}

Here we discuss a way of constructing from a distribution another one having a higher probability to observe outliers. We call this procedure "put tail down".

\vspace{0.2cm}
Let $F(x)$ be a probability distribution function of random variable $X$ having finite second moment $\sigma^2$ and such that $F(-x) = 1-F(x)$ for all $x \in \R^1$. Take a parameter $p \in (0,1)$ and fix it. Define a new function 
\[ F_p(x)=(1-p)F(x) + p H(x),\]
where $H(x) = 0$ for $x<0$, and $H(x)=1$ for $x>0$. It is clear that $F_p(x)$ is probability distribution function for any $p \in (0,1)$. Of course, $F_p$ also has finite second moment $\sigma_p^2$, and $F_p(-x)=1-F_p(x)$. However, $\sigma_p^2=(1-p)\sigma^2 ,\sigma^2$. Let $Y_p$ be a random variable with probability distribution function $F_p$. Then
\[ \p\{|Y_p|>k\sqrt{1-p}\sigma\} = 2 \p\{Y_p>k\sqrt{1-p}\sigma\}=\] \[=2(1-p)\bigl(1-F(k\sqrt{1-p}\sigma)\bigr).\]

Denoting $\bar{F}(x)=1-F(x)$ rewrite previous equality in the form
\begin{equation}
\label{eq3a}
\p\{|Y_p|>k\sqrt{1-p}\sigma\} = 2 (1-p) \bar{F}(k\sqrt{1-p}\sigma).
\end{equation}
For $Y_p$ to have more outliers than $X$ it is sufficient that
\begin{equation}\label{eq4a}
(1-p) \bar{F}(k\sqrt{1-p}\sigma) > \bar{F}(k\sigma).
\end{equation}
There are many cases in which inequality (\ref{eq4a}) is true for sufficiently large values of $k$. Let us mention two of them.
\begin{enumerate}
	\item Random variable $X$ has exponential tail. More precisely, 
	\[ \bar{F}(x) \sim C e^{-a x}, \;\text{as}\; x \to \infty , \]
	for some positive constants $C$ and $a$. In this case, inequality (\ref{eq4a}) is equivalent for sufficiently large $k$ to
	\[(1-p) > Exp\{-a\cdot k \cdot \sigma \cdot (1-\sqrt{1-p})\},\]
	which is obviously true for large $k$.
	\item $F$ has power tail, that is $\bar{F}(x) \sim C/x^{\alpha}$, where $\alpha>2$ in view of existence of finite second moment. Simple calculations show that (\ref{eq4a}) is equivalent as $k \to \infty $ to
	\[ (1-p)^{1-\alpha/2} <1. \]
\end{enumerate}

	The last inequality is true for $\alpha >2$. 

\vspace{0.2cm}
Let us note that the function $F_p$ has a jump at zero. However, one can obtain similar effect without such jump by using a smoothing procedure, that is by approximating $F_p$ by smooth functions.

"Put tail down" procedure allows us to obtain more outliers in view 
of two its elements. First element consists in changing the tail by smaller, but proportional to previous with coefficient $1-p$. The second element consist in moving a part of mass into origin (or into a small neighborhood of it), which reduces the variance.

\vspace{0.2cm} The procedure described above shows us that the presence of outliers may have no connection with existence of heavy tails of underlying distribution or with experimental errors. 

\section{Back to heavy tails. Estimation of tail index}\label{s3}
\setcounter{equation}{0}
As it has been mentioned above, in Section \ref{s1}, it is impossible to estimate tail index in general situation. However, it seems to be possible to construct upper (or lower) statistical estimators of tail index inside a special class of probability distributions. But {\bf what class of distributions allows such estimators?}

\vspace{0.2cm} To find such class let us consider a problem which seems (from the point of applications) to be far from the theory of heavy-tailed distributions. It appears in Medicine and considers a presence or absence of ``cure." 

\vspace{0.2cm} The probability of cure, variously referred to as the cure rate or the surviving fraction, is defined as an asymptotic value of the improper survival function  as time tends to infinity. 

Let X denote observed survival time. 
Statistical inference on
cure rates relies on the fact that any improper survival function $S(t) = \p\{X \geq t\}$ can be represented in the form:
\begin{equation}\label{eq5}
S(t) = a + (1-a)S_o(t),
\end{equation}
where $a = \p \{X=\infty \}$ is the probability of cure, and $S_o(t)$ is defined as the survival function for the time to failure conditional upon ultimate failure, i.e. 
\[ S_o(t) = \p\{ X \ge t | X<\infty \}. \]

\vspace{0.2cm}
Of course,
\[ a= \lim_{t \to \infty}S(t).\]
However, this relation cannot be used to construct any statistical estimator for the probability of cure. To have such a possibility we need to restrict the set of survival functions $S_o (t)$ under consideration to  {\bf a class of the functions with known speed (or known upper boundary of speed) of convergence to zero at infinity}. 

One of such classes is the set of distributions having increasing in average rate function (IFRA). More precisely, {\bf a distribution $F(x)$ concentrated on positive semi-axis belongs to the class IFRA if and only if the function
\[ -\frac{1}{x}\log (1-F(x)) \]
increases in $x \geq 0$} (see \cite{BP}).

\vspace{0.2cm} If $F$ belongs to the class IFRA then for any $t$ and $x$ such that $0<t \le x$
\[ -\frac{1}{x}\log (1-F(x)) \ge -\frac{1}{t}\log (1-F(t))  \]
that is
\begin{equation}\label{eq6}
1-F(x) \le \bigl(1-F(t) \bigr)^{x/t}.
\end{equation}
In other words, if we know the value of $F(t)$ then we have upper bound for the speed of convergence of $1-F(x)$ to zero as $x \to \infty$. This speed boundary (\ref{eq6}) is exponential.
 
Of course, one can construct statistical estimator for $F(t)$ using empirical distribution function. This allows one to obtain a lower bound for cure probability. However, our aim in this talk is not a study of cure, but the study of heavy tails. Therefore, we omit any estimators of cure probability, and go back to heavy-tailed distributions.

\vspace{0.2cm} To continue such study we need a modification of the hazard rate notion (see \cite{KYa}). Let $\varphi (u)$ be a nonnegative strictly monotonically decreasing function defined for all $u \ge 0$. Suppose in addition that its first derivative $\varphi^{\prime}$ is continuous and $\varphi (0) = -\varphi^{\prime}(0) =1$. We define the {\bf $\varphi$-hazard rate $r(t)=r_S(t)$ for the survival function $S(t)$ } by the following relations:
\[ \rho(t)=\rho_S(t)=\frac{d}{dt} \varphi^{-1}(S(t)), \]
\begin{equation}\label{eq7}
r(t)=r_S(t)=\rho (e^t)e^t .
\end{equation} 
We say, {\bf $F(t)$ belongs to the class $\varphi$-IFRA if and only if the function $r_S(t)$ increases in $t>0$, where $S(t)=1-F(t)$.}

\begin{thm}\label{th2} {\it \ref{th2}. Suppose that $X$ is a positive r.v. whose distribution function $F(x)$ belongs to the class $\varphi$-IFRA. Then for any $u>v>0$ holds
\begin{equation}\label{eq8}
S(u) \le \varphi \Bigl( \frac{\log u}{\log v} \varphi^{-1}(S(v))\Bigr),
\end{equation}
where $S(u)=1-F(u)$.}
\end{thm}

\vspace{0.2cm} Let us mention a particular case of Theorem \ref{th2}, when $\varphi(t)= \exp\{-t\}$. In this situation, class $\varphi$-IFRA coincides with the set of all distributions whose survival function $S(x)$ are such that $S(e^x)$ belongs to classical class IFRA. The inequality (\ref{eq8}) gives us
\[ S(u) \le u^{\log S(v)/\log v}, \; \; u>v>1. \]

\vspace{0.2cm}
Changing the restriction ``$r_S(t)$ increases in $t>0$'' in the definition of $\varphi$-IFRA class by ``$r_S(t)$ decreases in $t>0$'' we obtain the definition  of $\varphi$-DFRA class. For distributions from this class the inequality (\ref{eq8}) has to be changed by the opposite.

\section{Concluding Remarks}
\begin{enumerate}
\item We have seen that heavy-tailed distribution may appear as natural models in some problems of physics, technique and social sciences. Many of such models remain outside of this talk, say, the problems of rating of scientific publications.
\item Statistical inferences for such distributions must be model-oriented. There are no universal statistical procedures for the set of all heavy-tailed distributions.
\item The notion of outliers seems to be not defined mathematically. On intuitive level, outliers may not be indicators of the presence of large variance in data or that of experimental errors.
\item Additionally to previous item, the notion  of outliers may be defined in different ways for various model. The presence of such outliers cannot be considered as something negative in their nature. Outliers just reflect some specific properties of the studied process.                                                                                                                                                                                                                                                                                    \end{enumerate} 

\section*{Acknowledgment}
The work was partially supported by Grant GACR 16-03708S.

\end{document}